\documentstyle[amstex,amscd,12pt]{article}
\makeatletter \oddsidemargin5pt \evensidemargin5pt \textwidth 16cm
\date{}

\newtheorem{theorem}{Theorem}[section]

\newtheorem{problem}[theorem]{Problem}

\newcommand{\re}{{\Bbb R}}

\newcommand{\lo}{\longrightarrow}

\begin{document}

\title{A  property of $C_p[0,1]$}

\author{Michael  Levin}

\maketitle
\begin{abstract}
 We prove that for every finite dimensional compact metric space $X$ there is
 an open  continuous linear surjection from $C_p[0,1]$ onto
 $C_p(X)$. The proof makes use of embeddings introduced by
 Kolmogorov  and Sternfeld in connection with Hilbert's 13th
 problem.
\bigskip
\\
{\bf Math. Subj. Class.:} 54C35, 54F45.
\end{abstract}
\begin{section}{Introduction}\label{intro}

 All spaces are assumed to be separable metrizable and maps
 continuous. A compactum is a metric compact space. For a space $X$,
 $C_p(X)$ denotes the space of continuous real-valued functions 
 equipped with the topology of pointwise convergence.
   We refer to   the topology of   pointwise convergence
   as the $C_p$-topology.

 Let $X \subset Y_1 \times \dots \times Y_k$ be an
 embedding of a compactum $X$ into the product
 of compacta $Y_1,\dots, Y_k$. Define the space $Z$
 as
 the disjoint union of $Y_1, \dots, Y_k$ and define  the linear
 transformation $L : C(Z) \lo C(X)$
 by $L(g)(x)=g(y_1)+g(y_2)+\dots+g(y_k)$ for $g \in C(Z)$
 and $x=(y_1, \dots , y_k) \in X$ where $y_1, \dots, y_k$
  are the coordinates of $x$ in the product
 $Y_1 \times \dots \times Y_k$.  We will call $L$ the induced
 transformation  of the embedding
 $X \subset Y_1 \times \dots \times Y_k$. It is obvious that
 $L$ is continuous in  both the uniform topology and
 the $C_p$-topology on the function spaces.

 An embedding $X \subset Y_1 \times \dots \times Y_k$ is called
 basic if the induced
 transformation  $L$ is surjective. Note that, in general,
 a surjective  linear transformation  of
 the function spaces  on compacta which is  continuous in both
  the uniform topology  and the $C_p$-topology is not necessarily   open
 in the $C_p$-topology \cite{l-l-p}.
 It was shown in \cite{l-l-p} that the transformation induced by 
 a basic embedding in the product  of two spaces ($k=2$)
 is open in the $C_p$-topology and it is an open problem
 if the similar result holds for $k>2$. 
  In this paper we  will give a partial answer to this problem for
 two types of basic embeddings, namely,  Kolmogorov
 and  Sternfeld-type embeddings

 Sternfeld \cite{yaki-dendrites}
  constructed  for every  $n$-dimensional compactum $X$
 a  basic embedding  of $X$ into the product
 of $n+1$ one-dimensional compacta.
 We will call this embedding
 a Sternfeld-type embedding. Sternfeld-type embeddings
   are defined in Section
  \ref{sternfeld-type}.

 Adjusting Kolmogorov's solution of Hilbert's 13th problem given in
 Kolmogorov's famous superposition theorem \cite{kolmogorov}, Ostrand
 \cite{ostrand}
 defined for every $n$-dimensional compactum $X$ a basic embedding
 $X \subset [0,1]^{2n+1}$ which we will call a Kolmogorov-type
 embedding. Kolmogorov-type embeddings
  are described in Section
 \ref{kolmogorov-type}.

 In Sections \ref{sternfeld-type}
 and \ref{kolmogorov-type} we will prove the following theorems.
  \begin{theorem}
  \label{sternfeld}
  Let $X \subset Y_1\times \dots \times Y_{n+1}$ be a Sternfeld-type
  embedding of an $n$-dimensional compactum $X$ into the product
  of one-dimensional compacta $Y_1, \dots, Y_{n+1}$. Then the induced
  transformation  is  open in the $C_p$-topology.
  \end{theorem}

  \begin{theorem}
  \label{kolmogorov}
  Let $X \subset [0,1]^3$ be a Kolmogorov-type embedding of
  a $1$-dimensional compactum $X$ into the cube.
  Then the induced transformation is open in
  the $C_p$-topology.
  \end{theorem}

  Let $X$ be $n$-dimensional and compact. By Theorem
  \ref{sternfeld}  there is a $1$-dimensional compactum
  $Z$(=the disjoint union of $Y_1, \dots, Y_{n+1}$) for which
  $C_p(Z)$ admits an open continuous linear transformation
  onto $C_p(X)$.  Let
 $\hat{Z}$=the disjoint union of three copies of $[0,1]$.
  By Theorem \ref{kolmogorov} there is
  an open  continuous  linear transformation from $C_p(\hat{Z})$
  onto $C_p(Z)$. Embed $\hat{Z}$ into $[0,1]$
  and take the restriction transformation from
  $C_p[0,1]$ to $C_p(\hat{Z})$ which is obviously
  surjective, open and continuous. Thus we obtain the main result
  of the paper.

  \begin{theorem}
  \label{main}
  For every finite dimensional compactum $X$ there is
  an open continuous linear transformation from $C_p[0,1]$
  onto $C_p(X)$.
  \end{theorem}
  
  Note that $C_p [0,1]$ and $C_p(X)$ for a compactum $X$ are
  not isomorphic(=linearly homeomorphic) if $\dim X >1$ and 
  in  many cases when $\dim X=1$.
  In contrast to the uniform topology, the existence of 
  an isomorphism between $C_p (X)$ and $C_p(Y)$ for
  compacta $X$ and $Y$ implies a great deal of similarity
  between $X$ and $Y$, in particular, it implies that
  $\dim X =\dim Y$ \cite{pavlovskii}.

  Theorem \ref{main} generalizes some previous results by
  Leiderman, Pestov,  Morris and the author
  \cite{l-m-p}, \cite{l-l-p}. Open problems and
  related results are discussed in Section \ref{problems}.

    \end{section}

 \begin{section}{Sternfeld-type embeddings}\label{sternfeld-type}
 Let $X$ be compact and $n$-dimensional.
 Sternfeld \cite{yaki-dendrites}
 showed that there are a decomposition
 $X=A_1\cup \dots \cup A_{n+1}$ of $X$ into $0$-dimensional
 subsets $A_1, \dots, A_{n+1}$  and  $1$-dimensional
 compacta $Y_1, \dots, Y_{n+1}$ such that $X$ admits an embedding
 $X \subset Y_1 \times \dots \times Y_{n+1}$
 having the following property:
 $x=p_i^{-1}(p_i(x))$
 for every projection $p_i : X \lo Y_i$ and $x \in A_i$.
 We will call such an embedding of $X$ a Sternfeld-type embedding
 with respect to a decomposition $A_1, \dots, A_{n+1}$ of $X$.
 The spaces $Y_i$ can be chosen to be dendrites
 \cite{yaki-dendrites}. A different way  of constructing
 Sternfeld-type embeddings can be derived from \cite{lev-israel}.

 Sternfeld \cite{yaki-dendrites} proved that any Sternfeld-type embedding
 is basic. Sternfeld's proof is based on Borel measures and it is
 not clear at all if it can be applied  to prove
 Theorem \ref{sternfeld}. In this paper we use  another more
 constructive approach
 which is described in \ref{approx-sternfeld} and  which also shows
 that Sternfeld-type embeddings are basic.

 \begin{subsection}{An approximation procedure}
 \label{approx-sternfeld}
 Let $X \subset Y_1\times \dots \times Y_{n+1}$ be
 a Sternfeld-type embedding of  a $n$-dimensional compactum $X$ with
 respect to a decomposition
  $X=A_1 \cup \dots \cup A_{n+1}, \dim A_i=0$. Let us describe an
  approximation procedure showing that the embedding of $X$ is
  basic. The case $n=0$ is trivial. Assume that $n>0$. 

  Let  $f : X \lo \re$ be continuous  and  $c>0$ such
  that $\|f \| < c$. Fix $\epsilon >0$ which will be determined
  later and which will depend only on $\|f\|$,  $c$ and $n$.
  Take a disjoint  family ${\cal V}_i$
  of open subsets of $Y_i$ such that  ${\cal V}_i$ covers
  $p_i(A_i)$ and diam$f(p_i^{-1}(V)) < \epsilon$ for every
  $V \in {\cal V}_i$. For every $i$ choose a finite subfamily
  ${\cal U}_i \subset p^{-1}_i({\cal V}_i)$ such that
  ${\cal U}={\cal U}_1 \cup \dots \cup {\cal U}_{n+1}$ covers $X$
  and the elements of $\cal U$ are non-empty.
  By $V_U$, $U \in
  {\cal U}_i$ we denote the set of ${\cal V}_i$ such that
  $U=p_i^{-1}(V_U)$.
 For every $U \in {\cal U}$
  take a  non-empty subset $F_U \subset X$
  closed in $X$ such that for
  ${\cal F}=\{ F_U : U \in {\cal U}\}$
  covers $X$ and
  fix  a point $x_U \in F_U$.  For every $U \in {\cal U}_i$
   take
  a continuous function $\phi_U : Y_i \lo [0,1]$ such that
  $\phi_U (Y_i\setminus  V_U)=0$ and
   $\phi_U(p_i(F_U))=1$.
  Define $g'_i : Y_i \lo \re$ as $g'_i=\sum_{U\in {\cal U}_i}
 \frac{1}{n+1} f(x_U)\phi_U$. Clearly $\|g'_i\|\leq
 \frac{1}{n+1}\|f\|<\frac{1}{n+1}c$.
  Let us show that for every $x \in X$
  we have

  ${}_{}$ ${}_{}$${}_{}$ ($*$)  ${}_{}$  ${}_{}$ ${}_{}$
  ${}_{}$${}_{}$ ${}_{}$ $|f(x)-\sum_i g'_i(y_i)| < \frac{n}{n+1}c$
\\
 where
   $y_i=p_i(x) \in Y_i$ are the coordinates of $x$.

   Assume that $|f(x)| \leq \epsilon$.
   Then for every $U \in \cal U$ such that $x \in U$ we have
   $|f(x_U)|<2\epsilon$ and hence
   $|g'_i(y_i)|<\frac{2}{n+1}\epsilon$ for every $i$. Then
   $|f(x)-\sum_i g'_i(y_i)| < \epsilon + 2\epsilon$.
   Thus taking $\epsilon < \frac{n}{3(n+1)}c$ we get that  ($*$) holds.

  Assume that $f(x) > \epsilon$.
  Then for every $U \in \cal U$ such that $x \in U$ we have
   $0<f(x_U)< f(x) +\epsilon$ and hence
   $0<g'_i(y_i)<\frac{1}{n+1}(f(x) +\epsilon)$ for every $i$.
  Note that there is $F_U$, $U\in {\cal U}_j$, such that
  $x \in F_U$ and hence
  $\frac{1}{n+1}(f(x)
  -\epsilon)<g'_j(x)< \frac{1}{n+1}(f(x)  +\epsilon)$.
  Then  $0<f(x) -g'_j(x)<
    \frac{n}{n+1}(f(x)+\epsilon)$ and
    $0<\sum_{i\neq j} g'_i(x)<
    \frac{n}{n+1}(f(x)+\epsilon)$.
   Hence $|f(x)-\sum_i g'_i(y_i)|=|(f(x) -g'_j(x)) -\sum_{i\neq j} g'_i(x)|<
    \frac{n}{n+1}(f(x)+\epsilon)$.
   Thus taking $\epsilon < c -\|f\|$ we get that  ($*$) holds.
 In a similar way we check the case $f(x) <- \epsilon$ and get
   that in all the cases ($*$) holds.

Recall that by $Z$ we denote the disjoint union of $Y_1,\dots,
   Y_{n+1}$.
   Define $g' : Z\lo \re $ by $g'|_{Y_i}=g'_i$. We have that
   $\|g'\| < \frac{1}{n+1} c$
   and  $\|f -L(g')\|< \frac{n}{n+1}c$.
   Applying the described above procedure iteratively one can
   construct a sequence of maps $g^{(t)} : Z \lo \re$ such that
   $\|g^{(t)}\|< \frac{1}{n+1}(\frac{n}{n+1})^{t-1}c$
   and $\| f -L(\sum_{ s=1}^{ t} g^{(s)})\|  < (\frac{n}{n+1})^t c$.
   Then for $g=\sum_{s=1}^\infty g^{(s)}$ we have $f=L(g)$ and
   hence the embedding of $X$ is basic.

 \end{subsection}

 \begin{subsection}{Proof of Theorem \ref{sternfeld}}
 \label{proof-sternfeld}

Suppose that $X \subset Y_1\times \dots \times Y_{n+1}$ is
 a Sternfeld-type embedding of  an $n$-dimensional compactum $X$ with
 respect to a decomposition
  $X=A_1 \cup \dots \cup A_{n+1}, \dim A_i=0$.

 Let $Z' $ be a finite subset of $Z=$the disjoint union of
 $Y_i$'s. Denote $Y'_i=Y_i \cap Z'$,  $X'_i=p_i^{-1}(Y'_i)\cap A_i$
 and $X'=X'_1 \cup \dots\cup X'_{n+1}$.
 It is clear that each $X'_i$ is finite and therefore $X'$
 is finite as well.

 Take any map $f : X \lo \re$
 such that $f(x)=0$ for every $x \in X'$.
 We will show that there is a map $g : Z \lo \re $ such that
 $g(z)=0$ for every $z \in Z'$ and $L(g)=f$. This
 property together with the fact that $L$ is open 
 in the uniform topology implies that $L$
 is open in
 the $C_p$-topology at the zero-map  on $Z$
  and hence, by the linearity of
 $L$, $L$ is open in
 the $C_p$-topology everywhere and Theorem \ref{sternfeld}  follows.

 The construction of $g$ follows the procedure described in
 \ref{approx-sternfeld} with the following additional
 requirements. We  assume that no point of $Y'_i \setminus p_i(A_i)$ is covered
 by ${\cal V}_i$,  $p_i^{-1}(V)$ contains
 at most one point of $X'$ for every $V \in {\cal V}_i$
 and if $U \in {\cal U}$  contains a point  of $X'$ then
 this point is also contained in $F_U$ and $x_U=F_U \cap X'$.
 It is easy to see that under these assumptions we have
 that $g'_i(y)=0$ for every $y \in Y'_i$ and
 $g'_i(p_i(x))=0$ for every $x \in X'$. Thus $g'(z)=0$
 for every $z\in Z'$ and $f(x)-L(g')(x)=0$ for every
 $x \in X'$. Hence the approximation procedure can be
 repeated iteratively  and we can construct
 a map $g : Z \lo \re$ with the required properties.
  The theorem is proved.

 \end{subsection}
 \end{section}

 \begin{section}{Kolmogorov-type embeddings}
 \label{kolmogorov-type}
 Let $X$ be an $n$-dimensional compactum. A cover of $X$ is
 said to cover $X$ at least $n+1$ times if every point of $X$
 belongs to at least  $n+1$ elements of the cover.
 Generalizing Kolmogorov's paper
 \cite{kolmogorov}, Ostrand \cite{ostrand}
 showed  that there is a countable family ${\Omega}$ of closed finite covers
  of $X$ such that inf$\{$mesh${\cal F}: {\cal F} \in \Omega \}=0$,
 each ${\cal F}\in \Omega$ covers $X$ at least $n+1$ times and
 each ${\cal F} \in \Omega$ splits into the union
 ${\cal F}={\cal F}_1 \cup \dots \cup {\cal F}_{2n+1}$
 of $2n+1$  families of disjoint sets. 
 Such a family of  covers ${\Omega}$ with
 a fixed splitting
 ${\cal F}={\cal F}_1 \cup \dots \cup {\cal F}_{2n+1}$
 for each ${\cal F}\in \Omega$ will be called
 a Kolmogorov family of covers (Kolmogorov  constructed
 a Kolmogorov family of covers   for $X=[0,1]^n$ and this family
  can be easily transferred  to
 an arbitrary
 $n$-dimensional compactum using a $0$-dimensional map to $[0,1]^n$).
 
 In the case when  a Kolmogorov family contains a cover of mesh$=0$
 (that may happen only if $X$ is finite) we assume that this cover
 appears in the family infinitely many times. 
 We will also assume that a Kolmogorov
 family contains only finitely many covers of 
 mesh$>\epsilon$ for every $\epsilon>0$.  Thus we always assume that 
 a Kolmogorov family is infinite and any infinite subfamily of 
 a Kolmogorov family is Kolmogorov as well.

 A map from $X$ to $[0,1]$ is said to separate a family of disjoint sets in $X$  
 if  the images of the sets are disjoint  in $[0,1]$. 
 Let ${\Omega}$ be a Kolmogorov family of covers of $X$. An embedding 
 $X \subset [0,1]^{2n+1}$ is said to separate  a cover ${\cal F} \in {\Omega}$
 if the projection $p_i : X \lo [0,1]$ separates ${\cal F}_i$ for every $i$. 
  An embedding 
 $X \subset [0,1]^{2n+1}$
 will be called a Kolmogorov-type embedding
 with respect to $\Omega$ if for every $\epsilon >0$ 
 there is 
 ${\cal F} \in \Omega$ with mesh${\cal F} < \epsilon$  such that 
 the embedding of $X$  separates $\cal F$.
  Note that almost all
 embeddings of $X$ into $[0,1]^{2n+1}$ are of Kolmogorov-type with
 respect to $\Omega$, see \cite{yaki-kolmogorov}. 
 In the next subsection we present an approximation procedure
 which  can be derived from  Kolmogorov's paper \cite{kolmogorov} and 
 which shows   that Kolmogorov-type
 embeddings are basic. This fact  was  observed by 
 Ostrand \cite{ostrand}, see also \cite{yaki-kolmogorov}.

 Note that for a Kolmogorov-type embedding $X \subset [0,1]^{2n+1}$ with respect to $\Omega$
 we can replace $\Omega$ by  any infinite  subfamily of $\Omega$  and the embedding of $X$ 
 will remain  of  Kolmogorov-type with respect to the replaced $\Omega$. (Thus, in particular,
 we may assume that  a Kolmogorov-type embedding with respect to $\Omega$
 separates every cover in $\Omega$.)

 \begin{subsection}{ An approximation procedure}
 \label{approx-kolmogorov}
 Let $X \subset [0,1]^{2n+1}$ be a Kolmogorov-type embedding with
 respect to a family of covers $\Omega$. Here we describe an
  approximation procedure showing that the embedding of $X$ is
  basic. The case $n=0$ is trivial. Assume that $n>0$.
  
  Let  $f : X \lo \re$ be continuous  and  $c>0$ such
  that $\|f \| < c$. Fix $\epsilon >0$ which will be determined
    later and which will depend only on $\|f\|$,  $c$ and $n$.
       Clearly we may assume that each $\cal F \in \Omega$
 consists of non-empty sets. 
  Choose any cover ${\cal F} \in \Omega$ such that 
  mesh$f({\cal F})=\{ f(F): F \in {\cal F} \}< \epsilon$.
   For every non-empty $F  \in \cal F$ fix a point $x_F \in X$ such  that
   $f(x_F)$ is at distance$<\epsilon$ from $f(F)$.
   (For  showing  that $X \subset [0,1]^{2n+1}$ is basic is enough to
   choose $x_F$ in $F$, however in the proof of Theorem \ref{kolmogorov}
   we may need to choose $x_F$ outside $F$.)
   
   Define $g'_i : [0,1] \lo \re$, $1\leq i \leq 2n+1$, such that
   $g'_i(F)=\frac{1}{2n+1}f(x_F)$ for every non-empty $F \in {\cal F}_i$ 
   and $\| g'_i \| \leq \frac{1}{2n+1}\| f\| < \frac{1}{2n+1}c$.
    Let us show that for every $x \in X$
  we have

  ${}_{}$ ${}_{}$${}_{}$ ($*$)  ${}_{}$  ${}_{}$ ${}_{}$
  ${}_{}$${}_{}$ ${}_{}$ $|f(x)-\sum_i g'_i(y_i)| < \frac{2n}{2n+1}c$
\\where
   $y_i=p_i(x) \in [0,1]$ are the coordinates of $x$.
   
   Indeed,  recall that ${\cal F}$ covers $x$ at least $n+1$ times. Choose
   a set $I_+ \subset \{1,2,\dots,2n+1\}$ containing exactly $n+1$ indices
   such that $x$ is covered by ${\cal F}_i $ for  every $i \in I_+$ and
   denote $I_-=\{1,2,\dots,2n+1\} \setminus I_+$.
   For every $i \in I_+$  there is $F \in {\cal F}_i$
   containing $x$ and  hence
   $|\frac{1}{2n+1}f(x)-g'_i(y_i)|=\frac{1}{2n+1}|f(x)-f(x_F)| < \frac{1}{2n+1}2\epsilon$.
    Then $|f(x)-\sum_i g'_i(y_i)|=
   |\sum_{i \in I_+}(\frac{1}{2n+1}f(x)- g'_i(y_i)) + 
   \sum_{i\in I_-}(\frac{1}{2n+1}f(x) -g'_i(y_i))|
    < \frac{n+1}{2n+1}2\epsilon +  n(\frac{1}{2n+1}\| f\| +\frac{1}{2n+1}\|f\|)=
    \frac{2(n+1)}{2n+1}\epsilon  +\frac{2n}{2n+1}\|f\|$.
          Thus taking $\epsilon < \frac{n}{n+1}(c -\|f\|)$ we get that  ($*$) holds.

 Denote by $Z$  the disjoint union of $2n+1$ copies  $Y_i=[0,1]$,
  $1 \leq i \leq 2n+1$, of the interval $[0,1]$.
   Define $g' : Z\lo \re $ by $g'|_{Y_i}=g'_i$. We have that
   $\|g'\| < \frac{1}{2n+1} c$
   and  $\|f -L(g')\|< \frac{2n}{2n+1}c$
   where $L$ is the linear transformation induced by the embedding
   $X \subset [0,1]^{2n+1}$.
   Applying the described above procedure iteratively one can
   construct a sequence of maps $g^{(t)} : Z \lo \re$ such that
   $\|g^{(t)}\|< \frac{1}{2n+1}(\frac{2n}{2n+1})^{t-1}c$
   and $\| f -L(\sum_{ s=1}^{ t} g^{(s)})\|  < (\frac{2n}{2n+1})^t c$.
   Then for $g=\sum_{s=1}^\infty g^{(s)}$ we have $f=L(g)$ and
   hence the embedding of $X$ is basic.

 \end{subsection}

 \begin{subsection}{Embeddings of
 $1$-dimensional compacta}
 \label{kolmogorov-1-dim}
 Let $X$ be a one-dimensional compactum and  $X \subset [0,1]^3$
 a Kolmogorov-type embedding with respect to a Kolmogorov
 family $\Omega$ of covers of $X$.  Denote $Y_i=[0,1]$, $1 \leq i \leq 3$
 and $Z=$the disjoint union of $Y_i$'s.  Recall that by
 $p_i$ we denote the projection 
 $p_i : X \lo Y_i$.
 \\\\
 {\bf Reserved and free points}. We say that  $z \in Y_i \subset Z$ is 
 a reserved point of $Z$ with
 respect to $\Omega$ if for all but finitely many
 covers ${\cal F}$ in $\Omega$,
 ${\cal F}_{i}$ intersects  $p^{-1}_{i}(z)$, that is
 there is $F \in {\cal F}_{i}$ such that
 $F$ intersects  $p^{-1}_{i}(z)$. Note that, since
 $p_{i}$ separates   ${\cal F}_{i}$, at most one element of ${\cal F}_i$
 can intersect $p^{-1}_{i}(z)$.

 A point   $z \in  Z$  is said to be strongly
 reserved with respect to $\Omega$ if  ${\cal F}_i$ intersects
  $p^{-1}_{i}(z)$ for every ${\cal F} \in \Omega$ and  the collection
  $\{ F : F \cap p_i^{-1}(z)\neq \emptyset,
  F \in {\cal F}_i, {\cal F} \in {\Omega}\}$ converges to a point
  $x\in X$, that is every neighborhood of $x$ in $X$ contains all but finitely many
  elements of the collection. We will say that 
  the point  $x$  witnessing the   reservation of
  $z$ or say that $z$ is reserved by $x$.

  A point $z \in Y_i \subset  Z$ 
   which is not reserved with respect to $\Omega$
  is said to be free with respect to $\Omega$.
  A point $z \in Y_i \subset  Z$ is said to be fully free with respect
  to $\Omega$ if 
     ${\cal F}_i$ does not intersect $p_i^{-1}(z)$ for every
  ${\cal F} \in \Omega$.
  
 It is obvious that 
  if $z\in Z$ is  reserved (strongly reserved, 
  reserved by $x \in X$, fully free) with respect to
   $\Omega$   then $z$ remains to be reserved 
   (strongly reserved, reserved by $x$, fully free respectively)
  with respect to any infinite subfamily of $\Omega$.
  Note that if $z \in Z$ is reserved (free) with respect to $\Omega$
  then replacing $\Omega$ by its infinite subfamily we can get 
  that $z$ is strongly reserved (fully free) with respect
  to $\Omega$.
   Also note that, since each ${\cal F} \in \Omega$
  covers $X$ at least twice,  every point $x \in X$  has at least 
  two coordinates reserved and   if every  coordinate of $x$ is 
  either strongly reserved or fully free with respect  $\Omega$ 
    then at least two coordinates
  of $x$ are reserved by $x$.
  \\\\
  {\bf Chains.} 
  A chain $\chi$ of length $m$ with respect to $\Omega$ is a couple $\chi=(A,B)$ such that
  $A=\{z_0, z_1,\dots, z_{2m}\}$ is a sequence of $2m+1$ elements of  $Z$,
  $B=\{x_1, x_2, \dots x_m \}$ is a sequence if $m$  elements of $X$
  such that  $z_{2j-2}, z_{2j-1}$ and   $z_{2j}$ are  the coordinates 
  of $x_j$, 
   the points $\{z_0,\dots, z_{2m-1}\}$
  are strongly reserved with respect to $\Omega$,  $z_{2j-2}$ and $z_{2j-1}$ are
  reserved by $x_j$ and the point $z_{2m}$ is strongly reserved or fully
  free with respect to $\Omega$. Note the coordinates 
  $z_{2j-2}, z_{2j-1}$ and   $z_{2j}$ of $x_j$ do not necessarily go 
   in the order corresponding to the order of the coordinates
   $x_j=(y_1, y_2, y_3) $  of $x_j$ in the product of  $Y_1$, $Y_2$ and $Y_3$
   (for example it may happen that $z_{2j-2}=y_2$),
   the only thing that we assume is that 
   $\{z_{2j-2}, z_{2j-1}, z_{2j} \}= \{y_1, y_2, y_3 \}$ as subsets of $Z$.
     The chain $ \chi$ is of length $0$
   if $A$ contains only one point $z_0$ and $B=\emptyset$.
   The points $z_0$ and $z_{2m}$ 
   ($x_1$ and $x_m$ if $m>0$) are called the initial   and 
   the terminal $Z$-points ($X$-points respectively)
    of the chain $\chi$.  The sequences $A$ and $B$ are called  the $Z$-sequence 
    and the $X$-sequence of $\chi$ respectively.
      A chain $\chi'=(A', B')$ is said to be an extension of the chain $\chi$
   if the sequences $A'$ and $B'$  starts with $A$ and $B$ respectively.
   A chain $\chi'$ is said to be a continuation of the chain $\chi$
   if the terminal $Z$-point of $\chi$ is the initial $Z$-point of $\chi'$.
   If $\chi'=(A', B')$ is a continuation of length $m'$
   of the chain $\chi (A, B)$ then we can define
   the chain $\chi''=\chi + \chi'=(A'',B'')$ of length$=m+m'$
   by letting  $B''$ be the sequence  $B$    followed by the sequence $B'$ and 
   $A''$ be the sequence  $A$  followed by the sequence  $A'$  
   when the last element of $A$ is identified with
   the first one of $A'$. We will call $\chi$ and $\chi'$ the head of length $m$
   and the tail of length $m'$ respectively of the chain $\chi''$ and also
   write $\chi=\chi'' - \chi'$ and  $\chi'=\chi'' - \chi$.
 \\\\  
   {\bf Almost  free, periodic  and non-periodic points.} 
     If the terminal $Z$-point of the chain $\chi$ is 
     not free then $\chi$ can be extended in the following way.
    Define $x_{m+1}=$ the point of $X$ witnessing the reservation of $z_{2m}$,
    replace $\Omega$ by its infinite subfamily in order to get that
    every coordinate of $x_{m+1}$ is either strongly reserved or fully free
    and define $z_{2m+1}$ and $z_{2m+1}$  as the  two other coordinates of $x_{m+1}$
    (in addition to $z_{2m}$) such that $z_{2m+1}$ is reserved  by $x_{m+1}$
    (recall that $x_{m+1}$ has at least two coordinates reserved by $x_{m+1}$).
    Note that constructing chains we will constantly replace by $\Omega$ its
    infinite subfamily. Thus talking about two (or finitely many) chains built
    one after another we always refer to the smallest subfamily  obtained
    in the last replacement.

    Let $Z'$ be a finite subset of $Z$. Enumerate  the points of $Z'$ 
    and 
    for every $z \in Z'$ define  the chain $\chi(z,0)$ as the chain of length $0$
    with the initial Z-point $z$.  For every $m$  we will replace $\Omega$
    by its infinite subfamily  and construct  for every $z \in Z'$ a chain 
    $\chi(z, m)$  proceeding from $m$ to $m+1$ as  follows.
          Define $\chi(z, m+1)=\chi(z,m)$ if the terminal $Z$-point of  
   $ \chi(z,m)$ is free.  According to our enumeration of
   $Z'$ go over the points $z$ of $Z'$ such that
  the terminal $Z$-point of  $ \chi(z,m)$ is not free and 
  replacing (if necessary) each time $\Omega$
   by its infinite subfamily extend  $\chi(z, m)$ to a chain
   of  $\chi(z, m+1)$  as described above (by adding one element to
   the $X$-sequence and two elements to $Z$-sequence of $\chi(z, m)$).
   Thus the length of $\chi(z, m) \leq m$ and the length of $\chi(z,m)=m$ 
   if the terminal $Z$-point of $\chi(z,m)$ is not free.

   Let us call $z \in Z'$ almost free  if there is $m$ such that
  the terminal $Z$-point of  $\chi(z,m)$ is free. If $z\in Z'$ is not
  almost free then we will say that $z$ is 
  periodic  if there  is $m$
   such that the $X$-sequence of  $\chi(z,m)$ contains  two equal elements,
   and $z$ is 
    non-periodic otherwise (that is
   for every $m$ all the elements of the $X$-sequence of $\chi(z,m)$  are distinct).
   
    Define  the $X$-support ($Z$-support)  of 
   $z \in Z'$ as the subset of $X$ ($Z$)
   consisting of all the elements of the $X$-sequences
   ($Z$-sequences)  of $\chi(z,m)$ for all $m$.
   The $Z$-support of $z$ is the union of the coordinates 
   of the points of the $X$-support of $z$. It is obvious that
   the $X$-support and the $Z$-support of an almost free  point in $Z'$ are
   finite.
   
   Note that if for two chains $\chi$ and $\chi'$ of 
   length $m$ and $m'$ respectively we have that $x_j=x'_{j'}$ for
  the elements $x_j $,  $x'_{j'} $, $j<m$ and $j' < m'$ in the $X$-sequences 
    of $\chi$ and $\chi'$ respectively then
   for the  elements $x_{j+1}$ and $x'_{j'+1}$ following 
   $x_j$ and $x'_{j'}$ in the $X$-sequences 
   we also have $x_{j+1}=x'_{j'+1}$. Indeed, if $x_j=x_{j+1}$ then all the coordinates
   of $x_j$ are reserved by $x_j$ and therefore $x'_{j'+1}=x_j$. If
   $x_{j+1} \neq x_j$ then one of the coordinates  of $x_j$ is reserved 
   by $x_{j+1}$  and therefore $x'_{j'+1}=x_{j+1}$. 
   Also note that if  we assume
   in addition that  $j+1 < m$, $j'+1 < m'$ and
   the elements of  the $X$-sequence of $\chi$ are  
   distinct  then not only $x_{j+1}=x'_{j'+1}$
   but also 
     $z_{2j}=z'_{2j'}, z_{2j+1}=z'_{2j'+1}$ and 
   $z_{2j+2}=z'_{2j'+2}$ for
    the corresponding elements of the $Z$-sequences of $\chi$ and $\chi'$
    respectively. Indeed,  $z_{2j}, z_{2j+1}, z_{2j+2}$ are the coordinates of $x_{j+1}$
     which are characterized by  the properties:
      $z_{2j+2}$ is the only coordinate of $x_{j+1}$ not reserved by    $x_{j+1}$, 
   and  $z_{2j}$ is the only coordinate of $x_{j+1}$ which is also a coordinate
   of $x_{j}$ (since otherwise either $z_{2j+1}$ or $z_{2j+1}$ would be
   reserved by $x_j$ and by   $x_{j+1}$ or $x_{j+2}$ which  contradicts the assumption 
   that $x_j$, $x_{j+1}$ and $x_{j+2}$ are distinct).
       Then the required equalities 
    follow  because the similar characterization holds for
    the coordinates of $x'_{j'+1}$ and $x_{j}=x'_{j'}$
    $x_{j+1}=x'_{j'+1}$, $x_{j+2}=x'_{j'+2}$.
   
     Thus we get the following facts for the points in $Z'$. 
     The $X$-support and the $Z$-support of 
     a periodic point  are finite.  The $X$-support  
      of a non-periodic point  $z$ is infinite and  the elements of the $X$-sequence 
      of $\chi(z, m)$ are distinct  for every $m$. The $X$-supports of a non-periodic point 
      and a periodic point   are disjoint. 
      The $X$-supports of a non-periodic point $z $
      and an almost free point are disjoint because otherwise
      the $X$-sequence of $\chi(z,m)$ would contain 
      an element with a fully free coordinate for some $m$. 
      Since
       each element of the $Z$-sequence of $\chi(z,m)$ of a non-periodic 
       or periodic point $z$ 
        is reserved 
       by an element of the $X$-sequence of $\chi(z,m+1)$  we get that
                the elements of the $Z$-sequence of $\chi(z,m)$ of 
                a non-periodic point $z$ 
       are distinct for every $m$ and 
        the $Z$-supports of a non-periodic point  and a periodic point 
      are disjoint. Because every point of the $Z$-support of an almost free
      point in $Z'$ is either fully free or reserved by a point in its $X$-support
      we  get that 
        the $Z$-supports of a non-periodic point  and an almost free point 
      are disjoint.

      It also follows from what we said before that if for two non-periodic points
      $z_1$ and $z_2$ in $Z'$ the chains $\chi(z_1, m_1)$ and
      $\chi(z_2, m_2)$, $m_1, m_2>0$, have the same terminal $Z$-points then
            $\chi(z_1, m_1+k)-\chi(z_1, m_1)= \chi(z_2, m_2+k)-\chi(z_2, m_2)$
      for every $k$ (the tails of length $k$ of $ \chi(z_1, m_1+k)$
       and $ \chi(z_2, m_2+k)$ are equal).
              Let us say that two non-periodic points in $Z'$
   are equivalent if their $X$-supports intersect. Then
   for every non-periodic point $z \in Z'$ we can  find $m(z)>0$
   such  that if $z_1, z_2$ are equivalent non-periodic points in $Z'$
   we have that the chains $\chi(z_1, m(z_1))$ and $\chi(z_2, m(z_2))$
   have the same terminal $Z$-points. Denote by $Z'_{-}$ a collection 
   non-periodic points having 
   exactly  one representative in each  equivalence class.
    For every $z \in Z'_{-}$ denote 
   $\chi_{-}(z, k)=\chi(z, m(z)+k)-\chi(z, m(z))$ and call 
   $\chi_{-}(z, k)$ the reduced chain of length $k$ generated by $z$.
   Note that for every $k$ all the elements of both the $X$-sequence and $Z$-sequence
   of $\chi_{-}(z, k)$, $z \in Z'_{-}$
    are distinct and  for $z_1, z_2 \in Z'_{-}$,
   $z_1 \neq z_2$, we have that
    the elements of the  $X$-sequence and $Z$-sequence of  $\chi_{-}(z_1,k)$ are distinct
   from the elements of the $X$-sequence  and $Z$-sequence respectively of $\chi_{-}(z_2,k)$
   
   Define  $X' \subset X$ as the set consisting of
   the $X$-supports of almost free and periodic points of $Z'$ and
   the elements of the $X$-sequences of the chains $\chi(z, m(z))$ for the non-periodic
   points $z \in Z'$. Similarly define $Z'_{+} \subset Z$ as the set
   consisting of the $Z$-supports of almost free and periodic points of $Z'$
   and the elements of the $Z$-sequences of the chains $\chi(z, m(z))$
   for the non-periodic points $z \in Z'$. 
   
   Clearly $X'$ and $Z'_{+}$ are finite and $Z'_{-} \subset Z' \subset Z'_{+}$.
    Fix an integer $k$  and let $0 \leq j \leq k$.  Define $X^{(j)}\subset X$ as the union
  of $X'$ and the elements of the $X$-sequences of $\chi_{-}(z, j)$  for all $z \in Z'_{-}$
   and define  $ Z^{(j)} \subset Z$ as 
   the union
   of $Z'_{+}$ and the elements of the $Z$-sequences of $\chi_{-}(z, j)$ for all $z \in Z'_{-}$.
   
     It follows from our construction that 
    the elements of the $X$-sequence of
   $\chi_{-}(z,k)$,  $z \in Z'_{-}$, do not lie in $X^{(0)}=X'$,
    the initial $Z$-point of 
   the chain $\chi_{-}(z,k)$, $z \in Z'_{-}$,  is the only element of the $Z$-sequence
   of  $\chi_{-}(z,k)$ lying 
   in $Z^{(0)}=Z'_{+}$, the coordinates  of the points of $X^{(j)}$ are contained in $Z^{(j)}$
   and every point of $Z^{(j)}$ is either fully free or reserved by a point 
   in $X^{(j+1)}$.

 \end{subsection}

 \begin{subsection}{Proof of Theorem \ref{kolmogorov}}
 \label{proof-kolmogorov}
 Let $X \subset [0,1]^3$ be a Kolmogorov-type embedding
 with respect to $\Omega$.
 We will show that the induced transformation is open at the zero-map
 on $Z$ and this proves  Theorem \ref{kolmogorov}.
 Take a finite subset $Z'$ of $Z$. Following 
 \ref{kolmogorov-1-dim} replace $\Omega$ by its infinite subfamily
 and define the finite sets  $X'\subset X$ and $Z'_{-}\subset Z' \subset Z'_{+} \subset Z$.  
    
       Let  a map $\phi : X \lo \re$  be such that $\phi(x)=0$ for
  every $x \in X'$ and let $\delta >0$. We will construct
  a map $g : Z \lo \re$ such  that
  $\|\phi -L(g)\| < \delta$ and $g(z)=0$
  for every $z \in Z'$. This shows that $L$
  is open at the zero-map on $Z$
  because $L$ is open in the uniform topology and Theorem \ref{kolmogorov} follows.
  The case $\| \phi \|=0$ is trivial so we can assume that
  $\| \phi \|>0$.

  Fix an integer $k \geq 0$ which will be determined later
  and which will depend only
  on $\|\phi\|$ and $\delta$. Again  replace  $\Omega$ by its
  infinite subfamily and following \ref{kolmogorov-1-dim}
  construct the chains $\chi_{-}(z,k), z \in Z'_{-}$, and
  the sets $X^{(j)}$ and $Z^{(j)}$, $0 \leq j \leq k$.
  
  Define $h : Z^{(k)} \lo \re$
  such that $h(z)=0$ for every $z \in Z'_{+}$ and
  for every chain $\chi_{-}(z,k)=(A,B)$, $z \in Z'_{-}$, with
  $A=\{x_1, \dots, x_k \}$ and $B=\{z_0,\dots , z_{2k} \}$ define
  $h(z_{2j-2})=0$ and 
  $h(z_{2j-1})=\phi(x_{j})$
   in order to get
  $\phi(x_j)=h(z_{2j-2})+h(z_{2j-1})+h(z_{2j-2})$ for $1\leq j \leq k$.
  Extend $h$ over $Z$ such that $\| h \| \leq \| f \|$.
  Then $\| L(h) \| \leq 3\|\phi\|$ and $\phi(x)=L(h)(x)$  for every $x \in X^{(k)}$.
  
  Thus for 
   $f=\phi-L(h)$ we have 
  $\| f \| \leq 4 \|\phi \|$ and $f(x)=0$ for every $x \in X^{(k)}$.
  Now we will apply the approximation procedure \ref{approx-kolmogorov} for $f$
  and $c=5\|\phi \|$
  to construct the map $g' : Z \lo \re$. The approximation procedure
  will be applied  with the following
  additional requirements. We assume that the cover ${\cal F}\in \Omega$ is chosen
  so that each element of $\cal F$ contains at most one point of $X^{(k)}$
  and we assume that for every $1 \leq i \leq 3$ we have that
  the set $p_i(F)$ contains at most one point of $Z^{(k)}\cap Y_i$ for
   every  $F \in {\cal F}_i$, $g'_i(z)=0$ for every fully  free $z \in Z^{(k)}\cap Y_i$
  and, finally, if $p_i(F)$, $F \in {\cal F}_i$,
   contains a strongly reserved point of $z \in Z^{(k)}\cap Y_i$ then
   $F$ is so close to the point $x\in X$ 
   witnessing the reservation of $z$ that we can set
  $x_F=x$
  (recall that every point of $Z^{(k)}$ is either fully free or strongly reserved).
  One can easily verify that such a cover ${\cal F}\in \Omega$ satisfying 
  the requirements of \ref{approx-kolmogorov} along with  those we just listed can
  be chosen indeed. The conditions we impose on $\cal F$ imply that
  that if  a point  $z\in Y_i \cap Z^{(k)}$ is reserved by a point
  $x \in X$ then $g'_i(z)=\frac{1}{2n+1} f(x)$.
  Thus we get that $g'(z)=0$ if $z\in Z^{(k)}$ is fully free and 
  $g'(z)=0$ if $z\in Z^{(k)}$ is reserved by a point $x \in X$ 
  with $f(x)=0$.

   Recall that every point of $Z^{(k-1)}$ is either fully free or reserved by
  a point of $X^{(k)}$.
   Then, since  $f(X^{(k)})=0$ 
  (that is $f(x)=0$ for every $x \in X^{(k)}$)
   we get  that $g'(Z^{k-1)})=0$. Since
  the coordinates of the points in $X^{(k-1)}$ are contained in $Z^{(k-1)}$
  we also get that $L(g')(X^{(k-1)})=0$. Thus applying the approximation
  procedure iteratively we can construct the maps $g^{(t)} : Z \lo \re$, $1\leq  t \leq k$,
  such that $L(g^{(t)})(X^{(k-t)})=0$, $g^{(t)}(Z^{(k-t)})=0$, 
   $\|g^{(t)}\|< \frac{1}{2n+1}(\frac{2n}{2n+1})^{t-1}c$
   and $\| f -L(\sum_{ s=1}^{ t} g^{(s)})\|  < (\frac{2n}{2n+1})^t c$.
   Then for $g=h+\sum_{t=1}^k g^{(t)}$ we have that $g(Z^{(0)})=0$ and
   $\| \phi - L(g)\| < (\frac{2n}{2n+1})^k c$. Now assume that   $k$  
   is taken such that
   $(\frac{2n}{2n+1})^k c =(\frac{2n}{2n+1})^k (5\|\phi \|) < \delta$.
   Thus we get that $g(Z')=0$  and 
    $\| \phi - L(g)\| < \delta$, and  the theorem follows.

 \end{subsection}

 \end{section}

 \begin{section} {Problems}\label{problems}
 As we already mentioned in Section \ref{intro},
  Theorems \ref{sternfeld} and \ref{kolmogorov}  are partial positive
 solutions of the following  open  problem which was posed in \cite{l-l-p}.
 
 \begin{problem}
 \label{problem1}
 Let $X \subset Y_1 \times \dots \times Y_k$ be a basic embedding of 
 a compactum $X$ into the product of compacta $Y_1,\dots, Y_k$.
 Is    the  induced transformation always open in the $C_p$-topology?
 \end{problem}
 
 Problem \ref{problem1} in its full generality seems to be difficult, therefore
 it is justified to  discuss some  cases of this problem  related 
 to the types of basic embedding considered in this paper.
 
  It was already mentioned  in Section \ref{intro} that  Problem \ref{problem1}
 has the affirmative answer  if   $k=2$ \cite {l-l-p}. 
 The case $k=2$ and Theorem \ref{sternfeld}
 can be considered in the following generalizing context. 
 
 Let $X$ and $Y_1,\dots , Y_k$
 be compact,
 $X \subset Y_1 \times \dots \times Y_k$ and
 $p_i : X \lo Y_i$ the projections.  Define
 $S_i(X)=\{ x\in X: p_i^{-1}(p_i(x))=x\}$, $E_i(X)=X \setminus S_i(X)$,
   $E(X)=E_1(X)\cap \dots \cap E_k(X)$, $E^1(X)=E(X)$
   and by induction $E^{t}(X)=E(E^{t-1}(X))$.
    Let us call the embedding of 
 $X$ a Sternfeld embedding of general type if
 there is $t$ such that $E^{t}(X)=\emptyset$.
 Sternfeld showed that any basic  embedding  into
 the product of two spaces is of Sternfeld's general type and
 any Sternfeld embedding of general type is basic  \cite{yaki-hilbert}.
 Theorem \ref{sternfeld} admits a relatively simple generalization for
  embeddings with $E(X)=\emptyset$. 
 It would be interesting to answer
 
 \begin{problem}
 \label{problem2}
 Is the  induced transformation   of
  any  Sternfeld embedding of general type 
   open in the $C_p$-topology?
 \end{problem}
 
 Note that 
 not every basic embedding is of Sternfeld's general type (for example
 no embedding of a circle $S^1$ into $[0,1]^k$ is 
 of Sternfeld's general type).
 
 In connection to Theorem \ref{kolmogorov}
   it  seems very interesting to address the case
  of Kolmogorov-type embeddings of compacta of $\dim>1$.
 
 \begin{problem}
 \label{problem3}
 Is the induced transformation of any Kolmogorov-type embedding open
 in the $C_p$-topology?
 \end{problem}
 
 Note that the $[0,1]$ interval in Theorem \ref{main} cannot be replaced
 by a $0$-dimensional compactum \cite{l-l-p}. Also note Theorem \ref{main}
 does not hold if $X$ is not strongly countable dimensional \cite{l-m-p}.
  This suggests
 the following problem.
 
 \begin{problem}
 \label{problem4}
 Characterize compacta $X$  admitting  a linear open continuous transformation
 from $C_p[0,1]$ onto $C_p(X)$.
 \end{problem}
 
 Problem \ref{problem4} is also unsettled   for not necessarily open
 transformations \cite{l-m-p}.

 \end{section}

Department of Mathematics\\
Ben Gurion University of the Negev\\
P.O.B. 653\\
Be'er Sheva 84105, ISRAEL  \\
e-mail: mlevine@math.bgu.ac.il\\\\
\end{document}